\documentclass[12pt]{article}
\usepackage{latexsym, amssymb}

\DeclareMathAlphabet\gothic{U}{euf}{m}{n}

\makeatletter
\def\eqnarray{\stepcounter{equation}\let\@currentlabel=\theequation
\global\@eqnswtrue
\tabskip\@centering\let\\=\@eqncr
$$\halign to \displaywidth\bgroup\hfil\global\@eqcnt\z@
  $\displaystyle\tabskip\z@{##}$&\global\@eqcnt\@ne
  \hfil$\displaystyle{{}##{}}$\hfil
  &\global\@eqcnt\tw@ $\displaystyle{##}$\hfil
  \tabskip\@centering&\llap{##}\tabskip\z@\cr}

\def\endeqnarray{\@@eqncr\egroup
      \global\advance\c@equation\m@ne$$\global\@ignoretrue}

\def\@yeqncr{\@ifnextchar [{\@xeqncr}{\@xeqncr[5pt]}}
\makeatother

\begin{document}
\bibliographystyle{tom}

\newtheorem{lemma}{Lemma}[section]
\newtheorem{thm}[lemma]{Theorem}
\newtheorem{cor}[lemma]{Corollary}
\newtheorem{voorb}[lemma]{Example}
\newtheorem{rem}[lemma]{Remark}
\newtheorem{prop}[lemma]{Proposition}
\newtheorem{stat}[lemma]{{\hspace{-5pt}}}
\newtheorem{obs}[lemma]{Observation}
\newtheorem{defin}[lemma]{Definition}

\newenvironment{remarkn}{\begin{rem} \rm}{\end{rem}}
\newenvironment{exam}{\begin{voorb} \rm}{\end{voorb}}
\newenvironment{defn}{\begin{defin} \rm}{\end{defin}}
\newenvironment{obsn}{\begin{obs} \rm}{\end{obs}}

\newenvironment{emphit}{\begin{itemize} }{\end{itemize}}

\newcommand{\gota}{\gothic{a}}
\newcommand{\gotb}{\gothic{b}}
\newcommand{\gotc}{\gothic{c}}
\newcommand{\gote}{\gothic{e}}
\newcommand{\gotf}{\gothic{f}}
\newcommand{\gotg}{\gothic{g}}
\newcommand{\gothh}{\gothic{h}}
\newcommand{\gotk}{\gothic{k}}
\newcommand{\gotm}{\gothic{m}}
\newcommand{\gotn}{\gothic{n}}
\newcommand{\gotp}{\gothic{p}}
\newcommand{\gotq}{\gothic{q}}
\newcommand{\gotr}{\gothic{r}}
\newcommand{\gots}{\gothic{s}}
\newcommand{\gotu}{\gothic{u}}
\newcommand{\gotv}{\gothic{v}}
\newcommand{\gotw}{\gothic{w}}
\newcommand{\gotz}{\gothic{z}}
\newcommand{\gotA}{\gothic{A}}
\newcommand{\gotB}{\gothic{B}}
\newcommand{\gotG}{\gothic{G}}
\newcommand{\gotL}{\gothic{L}}
\newcommand{\gotS}{\gothic{S}}
\newcommand{\gotT}{\gothic{T}}

\newcommand{\mn}{\marginpar{\hspace{1cm}*} }
\newcommand{\mnn}{\marginpar{\hspace{1cm}**} }

\newcommand{\mnq}{\marginpar{\hspace{1cm}*???} }
\newcommand{\mnnq}{\marginpar{\hspace{1cm}**???} }

\newcounter{teller}
\renewcommand{\theteller}{\Roman{teller}}
\newenvironment{tabel}{\begin{list}%
{\rm \bf \Roman{teller}.\hfill}{\usecounter{teller} \leftmargin=1.1cm
\labelwidth=1.1cm \labelsep=0cm \parsep=0cm}
                      }{\end{list}}

\newcounter{tellerr}
\renewcommand{\thetellerr}{(\roman{tellerr})}
\newenvironment{subtabel}{\begin{list}%
{\rm  (\roman{tellerr})\hfill}{\usecounter{tellerr} \leftmargin=1.1cm
\labelwidth=1.1cm \labelsep=0cm \parsep=0cm}
                         }{\end{list}}
\newenvironment{ssubtabel}{\begin{list}%
{\rm  (\roman{tellerr})\hfill}{\usecounter{tellerr} \leftmargin=1.1cm
\labelwidth=1.1cm \labelsep=0cm \parsep=0cm \topsep=1.5mm}
                         }{\end{list}}

\newcommand{\Ni}{{\bf N}}
\newcommand{\Ri}{{\bf R}}
\newcommand{\Ci}{{\bf C}}
\newcommand{\Ti}{{\bf T}}
\newcommand{\Zi}{{\bf Z}}
\newcommand{\Fi}{{\bf F}}

\newcommand{\proof}{\mbox{\bf Proof} \hspace{5pt}} 
\newcommand{\remark}{\mbox{\bf Remark} \hspace{5pt}}
\newcommand{\ruimte}{\vskip10.0pt plus 4.0pt minus 6.0pt}

\newcommand{\simh}{{\stackrel{{\rm cap}}{\sim}}}
\newcommand{\ad}{{\mathop{\rm ad}}}
\newcommand{\Ad}{{\mathop{\rm Ad}}}
\newcommand{\Aut}{\mathop{\rm Aut}}
\newcommand{\arccot}{\mathop{\rm arccot}}
\newcommand{\capp}{{\mathop{\rm cap}}}
\newcommand{\rcapp}{{\mathop{\rm rcap}}}
\newcommand{\Capp}{{\mathop{\rm Cap}}}
\newcommand{\diam}{\mathop{\rm diam}}
\newcommand{\divv}{\mathop{\rm div}}
\newcommand{\dist}{\mathop{\rm dist}}
\newcommand{\codim}{\mathop{\rm codim}}
\newcommand{\RRe}{\mathop{\rm Re}}
\newcommand{\IIm}{\mathop{\rm Im}}
\newcommand{\Tr}{{\mathop{\rm Tr}}}
\newcommand{\Vol}{{\mathop{\rm Vol}}}
\newcommand{\card}{{\mathop{\rm card}}}
\newcommand{\supp}{\mathop{\rm supp}}
\newcommand{\sgn}{\mathop{\rm sgn}}
\newcommand{\essinf}{\mathop{\rm ess\,inf}}
\newcommand{\esssup}{\mathop{\rm ess\,sup}}
\newcommand{\Int}{\mathop{\rm Int}}
\newcommand{\Leibniz}{\mathop{\rm Leibniz}}
\newcommand{\lcm}{\mathop{\rm lcm}}
\newcommand{\loc}{{\rm loc}}

\newcommand{\mod}{\mathop{\rm mod}}
\newcommand{\spann}{\mathop{\rm span}}
\newcommand{\one}{1\hspace{-4.5pt}1}

\newcommand{\DWR}{}

\hyphenation{groups}
\hyphenation{unitary}

\newcommand{\tfrac}[2]{{\textstyle \frac{#1}{#2}}}

\newcommand{\cb}{{\cal B}}
\newcommand{\cc}{{\cal C}}
\newcommand{\cd}{{\cal D}}
\newcommand{\ce}{{\cal E}}
\newcommand{\cf}{{\cal F}}
\newcommand{\ch}{{\cal H}}
\newcommand{\ci}{{\cal I}}
\newcommand{\ck}{{\cal K}}
\newcommand{\cl}{{\cal L}}
\newcommand{\cm}{{\cal M}}
\newcommand{\cn}{{\cal N}}
\newcommand{\co}{{\cal O}}
\newcommand{\cs}{{\cal S}}
\newcommand{\ct}{{\cal T}}
\newcommand{\cx}{{\cal X}}
\newcommand{\cy}{{\cal Y}}
\newcommand{\cz}{{\cal Z}}

\newcommand{\wtozp}{W^{1,2}\raisebox{10pt}[0pt][0pt]{\makebox[0pt]{\hspace{-34pt}$\scriptstyle\circ$}}}
\newlength{\hightcharacter}
\newlength{\widthcharacter}
\newcommand{\covsup}[1]{\settowidth{\widthcharacter}{$#1$}\addtolength{\widthcharacter}{-0.15em}\settoheight{\hightcharacter}{$#1$}\addtolength{\hightcharacter}{0.1ex}#1\raisebox{\hightcharacter}[0pt][0pt]{\makebox[0pt]{\hspace{-\widthcharacter}$\scriptstyle\circ$}}}
\newcommand{\cov}[1]{\settowidth{\widthcharacter}{$#1$}\addtolength{\widthcharacter}{-0.15em}\settoheight{\hightcharacter}{$#1$}\addtolength{\hightcharacter}{0.1ex}#1\raisebox{\hightcharacter}{\makebox[0pt]{\hspace{-\widthcharacter}$\scriptstyle\circ$}}}
\newcommand{\scov}[1]{\settowidth{\widthcharacter}{$#1$}\addtolength{\widthcharacter}{-0.15em}\settoheight{\hightcharacter}{$#1$}\addtolength{\hightcharacter}{0.1ex}#1\raisebox{0.7\hightcharacter}{\makebox[0pt]{\hspace{-\widthcharacter}$\scriptstyle\circ$}}}

\newpage

 \thispagestyle{empty}
 
 \begin{center}
 \vspace*{-1.0cm}
\vspace*{1.5cm}
 
{\Large{\bf Hardy inequalities, Rellich inequalities}  }\\[3mm] 
{\Large{\bf and local Dirichlet forms  }}  \\[4mm]
\large Derek W. Robinson$^\dag$ \\[1mm]

\normalsize{January 2017}
\end{center}

\vspace{+5mm}

\begin{list}{}{\leftmargin=1.7cm \rightmargin=1.7cm \listparindent=15mm 
   \parsep=0pt}
   \item
{\bf Abstract} $\;$ 
First the Hardy and Rellich inequalities are defined for the submarkovian operator
associated with a local Dirichlet form.
Secondly,  two general conditions are derived which are sufficient  to deduce the Rellich inequality from the Hardy 
inequality.
In addition the Rellich constant is calculated from the Hardy constant.
Thirdly, we establish that the criteria for the Rellich inequality are verified for a large class of weighted second-order operators  on a domain $\Omega\subseteq \Ri^d$.
The  weighting near the boundary $\partial \Omega$  can be different from the weighting at infinity.
Finally these results are applied to weighted second-order operators on $\Ri^d\backslash\{0\}$ and to a general class of operators
of Grushin type.

\end{list}

\vfill


\vspace{0.5cm}

\noindent
\begin{tabular}{@{}cl@{\hspace{10mm}}cl}
$ {}^\dag\hspace{-5mm}$&   Mathematical Sciences Institute (CMA)    &  {} &{}\\
  &Australian National University& & {}\\
&Canberra, ACT 0200 && {} \\
  & Australia && {} \\
  &derek.robinson@anu.edu.au
 & &{}\\
\end{tabular}

\newpage

\setcounter{page}{1}

\newpage

\section{Introduction}\label{S1}

Our intent is twofold. 
First we  analyze Hardy and Rellich inequalities in the general framework of  local Dirichlet forms.
Secondly we apply the analysis to a large class of  divergence-form elliptic operators on domains of $\Ri^d$.
Our principal results give  verifiable criteria that allow the deduction of the Rellich inequality from the 
Hardy inequality.

There is an enormous literature concerning variants of the Hardy inequality and their applications
but somewhat less for the Rellich inequality.
We refer to the book by Balinsky, Evans and Lewis \cite{BEL} and the thesis by Ward \cite{War}  for 
background information and an indication of the relevant literature.
In order to explain our results we first establish some notation and recall some basic elements of the theory of Dirichlet forms.
We mostly adopt the definitions and terminology of Bouleau and Hirsch \cite{BH} (see also \cite{FOT}).
Subsequently we turn to the examination of elliptic operators on domains of Euclidean space.

Let $X$ denote  a locally compact $\sigma$-compact
metric space and $\mu$ a positive Radon measure with $\supp \mu = X$.
The corresponding real $L_p$-spaces are denoted  by $L_p(X)$.
Let $\ce$ denote a Dirichlet form with domain $D(\ce)$ on $L_2(X)$ and 
$H$  the self-adjoint submarkovian operator canonically associated with $\ce$.
Set  $B(\ce)=D(\ce)\cap L_\infty(X)$.
Then $B(\ce)$ is an algebra and a core of $D(\ce)$. 
Further let $B_{\rm loc}(\ce)$ denote the corresponding space of functions which are locally in $B(\ce)$, i.e.\ the space of $\mu$-measurable functions $\psi$
such that for every compact subset $K$ of $X$ there is a $\hat \psi \in B(\ce)$
with $\psi|_K = \hat \psi|_K$.
Next let $B_c(\ce)$ denote the subalgebra of $B(\ce)$ formed by the functions with compact support and set $C_c(\ce)=B_c(\ce)\cap C(X)$.
We assume that $C_c(\ce)$ is dense in $C_0(X)$, the space of continuous functions over $X$ which vanish at infinity,  with respect to the supremum norm   and  that  it is also  dense in $B_c(\ce)$ with respect to the $D(\ce)$-graph norm, i.e.\ the norm $\varphi\in D(\ce)\mapsto
\|\varphi\|_{D(\ce)}=(\ce(\varphi)+\|\varphi\|_2^2)^{1/2}$.
In addition  we assume  that $\ce$ is local in the sense of Bouleau and Hirsch \cite{BH}, Section~I.1.5.
In particular $\ce$ is local if $\ce(\varphi, \psi)=0$ for all $\varphi, \psi\in D(\ce)$ for which there is an $a\in\Ri$ such that 
$(\varphi+a\one)\,\psi=0$. 
(A slightly more specific property is introduced in \cite{FOT}, Section~1.1 and is referred to as strong locality.)
Finally, for each positive $\xi\in B(\ce)$ we define the truncated form $\ce_\xi$ by $D(\ce_\xi)=B(\ce)$  and
\begin{equation}
\ce_\xi(\varphi)=\ce(\varphi, \xi\varphi)-2^{-1}\ce(\xi,\varphi^2)
\label{eloc1.3}
\end{equation}
for all $\varphi\in B(\ce)$.
The truncated forms satisfy the Markovian properties characteristic of Dirichlet forms but are not necessarily closed.
Moreover, $\xi\mapsto \ce_\xi(\varphi)$ is a positive linear functional for each $\varphi\in B(\ce)$
and $\ce_\xi(\varphi)\leq \|\xi\|_\infty\,\ce(\varphi)$ for all $\varphi\in B(\ce)$ (see \cite{BH}, Proposition~4.1.1 for these properties).
Consequently, for each $\varphi\in B(\ce)$, the function $\xi\mapsto \ce_\xi(\varphi)$ extends by continuity to $C_0(X)$.
Then there is a positive Radon measure $\mu_\varphi$, the energy measure, such that $\mu_\varphi(\xi)=\ce_\xi(\varphi)$
for all $\xi\in C_0(X)$.
Note that if $\xi\in B(\ce)$ has compact support  one can also define $\ce_\xi$ on $B_{\rm loc}(\ce)\cap L_\infty(X)$ by setting $\ce_\xi(\varphi)=\ce_\xi(\hat\varphi)$
where $\hat\varphi\in B(\ce)$ is such that $\hat\varphi|_{\supp\xi}=\varphi|_{\supp\xi}$.
The definition is consistent by locality.

\smallskip

Next  let $\eta\in B_{\rm loc}(\ce)$. 
Then the  Dirichlet form  $\ce$ is defined to satisfy the $\eta$-Hardy inequality if
$\eta \,D(\ce)\subseteq L_2(X)$ and 
\begin{equation}
\ce(\varphi)\geq (\eta\,\varphi, \eta\,\varphi)
\label{eloc1.1}
\end{equation}
for all $\varphi\in D(\ce)$.
Since this condition is invariant under the map $\eta\mapsto |\eta|$ one may always assume that $\eta$ is positive.
Similarly,  $\ce$ is defined to satisfy the $\eta$-Rellich inequality if $\eta^2D(H)\subseteq L_2(X)$ and there is a $\sigma>0$ such that 
\begin{equation}
(H\varphi, H\varphi)\geq \sigma\,(\eta^2\varphi, \eta^2\varphi)
\label{eloc1.2}
\end{equation}
for all $\varphi\in D(H)$.

Our main result, which is proved in Section~\ref{S2},  establishes conditions which ensure that the $\eta$-Rellich inequality follows from the $\eta$-Hardy inequality. 

\begin{thm}\label{tloc1.1}
Assume $\ce$ is a local Dirichlet form and $\eta\in B_{\rm loc}(\ce)$  is positive.
Further assume
\begin{tabel}
\item\label{tloc1.1-1}
$\ce$ satisfies the $\eta$-Hardy inequality $(\ref{eloc1.1})$,
\item\label{tloc1.1-2}
there is a $\gamma\in \langle0,1\rangle$ such that $\ce_{\varphi^2}(\eta)\leq \gamma\,(\eta^2\varphi, \eta^2\varphi)$
for all $\varphi\in B_c(\ce)$,
\item\label{tloc1.1-3}
there exists a net $\{\rho_\alpha\}$ with $\rho_\alpha\in B_c(\ce)$ such that $0\leq \rho_\alpha\leq 1$,
\begin{equation}
\hspace{-2cm}\lim_\alpha(\varphi, \rho_\alpha\varphi)=(\varphi, \varphi) \;\;\;\;\;\;{ and}\;\;\;\;\; \lim_\alpha\ce_{\varphi^2}(\rho_\alpha)=0
\label{eloc1.4}
\end{equation}
for all $\varphi\in B(\ce)$.
\end{tabel}

It then follows that $\ce$ satisfies the $\eta$-Rellich inequality $(\ref{eloc1.2})$ with $\sigma=(1-\gamma)^2$.
\end{thm}

The Hardy inequality (\ref{eloc1.1}) is the quadratic form  expression of the ordering  $H\geq \eta^2$ of the self-adjoint operators
$H$ and  $\eta^2$ where the latter is interpreted as a multiplication operator.
Similarly, the Rellich inequality (\ref{eloc1.2}) corresponds to the order relation $H^2\geq\sigma\,\eta^4$.
Note that  the order relation is, however,  not  generally respected by taking squares unless the operators commute.
This is the role played by Condition~\ref{tloc1.1-2}; it imposes restrictions on the   commutativity of $H$ and $\eta$.
The condition can be rephrased in terms of the  energy  measures $\mu_\varphi$ and it is most transparent if 
these measures are absolutely continuous with respect to $\mu$.
The corresponding Radon--Nikodym derivatives $\Gamma\colon\varphi\in B(\ce)\mapsto\Gamma(\varphi)\in L_1(\Omega\,;\mu)$, which  are usually referred to as the {\it carr\'e du champ}, define  a  positive quadratic form  whose basic properties are developed in \cite{BH} Section~I.1.4.
Then one has
\[
\ce_{\varphi^2}(\eta)=\int_X d\mu\,\Gamma(\eta)\,\varphi^2=(\varphi, \Gamma(\eta)\varphi)
\]
for all $\varphi\in B_c(\ce)$.
Thus Condition~\ref{tloc1.1-2} is equivalent to the bounds
\[
0\leq \Gamma(\eta)\leq \gamma\,\eta^4
\;.
\]
But in applications to elliptic operators $\Gamma(\eta)$ is a measure of non-commutation.
For example, the {\it carr\'e du champ} corresponding to  the Laplacian $\Delta$ on $\Ri^d$ 
is given by $\Gamma(\eta)=|\nabla\eta|^2$ and formally $|\nabla\eta|^2=-2^{-1}\,[\, [\,\Delta,\eta\,],\eta\,]$.
Thus in this  case Condition~\ref{tloc1.1-2} leads to the bounds
\[
0\leq -2^{-1}\,[\, [\,\Delta,\eta\,],\eta\,]\leq \gamma\,\eta^4
\]
 on the double commutator of  $\Delta$ and $\eta$.
This  restriction on the commutation is the essential  content of Condition~\ref{tloc1.1-2} of the theorem.
Double commutator estimates of this type occur in many disparate areas of  mathematical physics and analysis, e.g.\ in quantum field theory,
\cite{GJ1} \cite{GJ} Section~19.4, \cite{RS2} Section~X.5, operator theory, \cite{Far} Section~II.12,   \cite{DrS},
elliptic equations \cite{Agm1}, elliptic regularity \cite{FP} \cite{ER27} \cite{RSi4},  etc.

 Condition~\ref{tloc1.1-3} of the theorem, which is independent of the Hardy--Rellich function $\eta$, corresponds to the existence of a special form of approximate identity $\{\rho_\alpha\}$  on $L_2(X)$.
Although it is not evident that an approximation of this type exists in general 
we do establish that it exists
for a large class of divergence-form elliptic operators  on a general domain  of $ \Ri^d$ if the operators satisfy an appropriate 
Hardy inequality.
To describe our results in the latter context we need some additional terminology.

Let $\Omega$ be a domain in $\Ri^d$, i.e.\ a connected open subset,
with boundary $\partial\Omega$ and equipped with the Euclidean metric 
$d(\,\cdot\,;\,\cdot\,)$.
Further let $x\in \Omega\mapsto d_\Omega(x)\in\langle0,\infty\rangle$ denote
the Euclidean distance to the boundary, i.e.\ $d_\Omega(x)=\inf_{y\in\Omega^c}d(x\,;y)$.
If  $c$ is  a strictly positive function on the half line $\langle0,\infty\rangle$ we define
$c_\Omega$ by $c_\Omega=c\circ d_\Omega$.
Then we define a Dirichlet form $h$ on $L_2(\Omega)$ as the closure of the form
\begin{equation}
\varphi\in C_c^\infty(\Omega)\mapsto h(\varphi)= \sum^d_{k=1} (\partial_k\varphi, c_\Omega\,\partial_k\varphi)
\;.
\label{eloc1.5}
\end{equation}
The form is closable, because $c$ is strictly positive,
and the closed form is automatically a local Dirichlet form  (see, for example, \cite{MR} Section~II.2.b).
Moreover, the form has a {\it carr\'e du champ} $\Gamma$ given by $\Gamma(\varphi)=c_\Omega\,|\nabla\varphi|^2$.
The   submarkovian operator $H$ corresponding to the form $h$ can be interpreted as the elliptic operator 
$-\sum^d_{k=1}\partial_k\,c_\Omega\,\partial_k$ with Dirichlet boundary conditions. 
In the context of the Hardy--Rellich inequality the choice $c(s)=s^\delta$ is conventional
but we will consider a broader class of coefficients and operators.

Our second principal result  is essentially  a corollary of Theorem~\ref{tloc1.1}.

\begin{thm}\label{tloc1.2}
 Let $c(s)=s^\delta\,(1+s)^{\delta'-\delta}$ with $\delta,\delta'\geq0$ and set
\[
\nu=\sup\{ |1-t/2|^2: \delta\wedge\delta'\leq t\leq \delta\vee\delta'\}
\;.
\]
Assume that the Dirichlet form $(\ref{eloc1.5})$  on $L_2(\Omega)$ corresponding to $c$ satisfies the Hardy inequality
\begin{equation}
h(\varphi)\geq a_1 \,(c_\Omega^{\,1/2}d_\Omega^{\;-1}\varphi, c_\Omega^{\,1/2}d_\Omega^{\;-1}\varphi)
\label{eloc1.6}
\end{equation}
for all $\varphi\in D(h)$ with $a_1>0$.

If $\nu<a_1$ then $H$ satisfies the Rellich inequality
\begin{equation}
(H\varphi, H\varphi)\geq a_2\, (c_\Omega \,d_\Omega^{\;-2}\varphi,c_\Omega \,d_\Omega^{\;-2}\varphi)
\label{eloc1.7}
\end{equation}
for all $\varphi\in D(H)$ with $a_2=(a_1-\nu)^2$.
\end{thm}

Our choice of the weight function $c$ in Theorem~\ref{tloc1.2} is dictated by its asymptotic behaviour.
The parameters $\delta$ and $\delta'$ govern the growth properties of $c_\Omega$ near the boundary and at infinity, respectively.
In particular one has  $\lim_{s\to0}c(s)\,s^{-\delta}=1$  and $\lim_{s\to\infty}c(s)\, s^{-\delta'}=1$.

Note that Condition~\ref{tloc1.1-1} of  Theorem~\ref{tloc1.1} is satisfied with $\eta^2=a_1\,c_\Omega\,d_\Omega^{\;-2}$ by 
the assumption that $h$ satisfies the Hardy inequality (\ref{eloc1.6}).
Moreover Condition~\ref{tloc1.1-2} of the earlier theorem  is not difficult to verify by direct calculation using the properties of $c$.
But the verification of the third condition of Theorem~\ref{tloc1.1} is more difficult.
We achieve this by adaptation of an argument introduced by Agmon \cite{Agm1} in his analysis of the exponential decay of solutions
of second-order elliptic equations. 
Agmon's arguments have earlier been used by Grillo \cite{Gril} to discuss Hardy and Rellich inequalities for  operators
constructed as sums of squares of vector fields.

Finally we note that the conclusions of the theorems are 
established for all functions  in the domain $D(H)$ of the submarkovian operator $H$.
Many derivations of the Rellich inequality are only valid on a core $D$ of the corresponding form $\ce$ but  not for a core of $H$.
In particular if $H_0$ is a symmetric elliptic operator on  a domain $\Omega\subseteq \Ri^d$ it is commonplace to use 
$D=C_c^\infty(\Omega)$ (see, for example, \cite{BEL}, Chapter~6).
Then the Dirichlet  form $\ce$ obtained by closure of $\varphi\in C_c^\infty(\Omega)\mapsto (\varphi, H_0\varphi)$ determines  the 
self-adjoint Friedrichs extension $H_{\!F}$ of $H_0$ but the closure of the form  $\varphi\in C_c^\infty(\Omega)\mapsto (H_0\varphi, H_0\varphi)$ determines the Friedrichs extension $(H_0^{\,2})_{F}$ of $H_0^{\,2}$ which  usually differs from $H_{\!F}^{\,2}$.
In general one has $(H_0^{\,2})_F\geq H_{\!F}^{\,2}$ with equality if and only if  $H_0$ is essentially self-adjoint.

\section{Locality estimates}\label{S2}

In this section we give the  proof of Theorem~\ref{tloc1.1}.
It is   based on several identities and estimates which are a direct result of the 
structure of the Dirichlet form $\ce$ and the locality condition.
We begin by collecting some specific results of relevance to the proof.

The locality property can be  exploited by use of Anderssen's representation theorem \cite{And} (see also \cite{Roth})
which is reformulated in \cite{BH}, Theorem~I.5.2.1. 
We will use the formulation given in \cite{ERSZ2} (see also \cite{AH}).

\begin{prop}\label{ploc2.1}
Let $\ce$ be a  local Dirichlet form on $L_2(X)$
and $\varphi_1,\ldots,\varphi_n \in B(\ce)$.
Then  there exists a unique real Radon measure 
$\sigma_{ij}^{(\varphi_1,\ldots,\varphi_n)}$
on $\Ri^n$ 
such that $\sigma_{ij}^{(\varphi_1,\ldots,\varphi_n)} = \sigma_{ji}^{(\varphi_1,\ldots,\varphi_n)}$
for all $i,j \in \{ 1,\ldots,n \} $, $\sum^n_{i,j=1}\xi_i\,\xi_j\,\sigma_{ij}\geq0$ for all $\xi_i,\xi_j\in\Ri$  and 
\begin{equation}
\ce(F_0(\varphi_1,\ldots,\varphi_n), G_0(\varphi_1,\ldots,\varphi_n))
= \sum_{i,j=1}^n \int_{\Ri^n} d\sigma_{ij}^{(\varphi_1,\ldots,\varphi_n)} \, 
     { \frac{\partial F}{\partial x_i} } \, \frac{\partial G}{\partial x_i}
\label{eloc2.1}
\end{equation}
for all $F,G \in C^1(\Ri^n)$ where $F_0=F-F(0)$ and $G_0=G-G(0)$.
Let $K$ be a compact subset of $\Ri^n$ such that 
$(\varphi_1(x),\ldots,\varphi_n(x)) \in K$ for almost every $x \in X$.
Then $\supp \sigma_{ij}^{(\varphi_1,\ldots,\varphi_n)} \subseteq K$ for all 
$i,j \in \{ 1,\ldots,n \} $.
In particular, if $i \in \{ 1,\ldots,n \} $ then $\sigma_{ii}^{(\varphi_1,\ldots,\varphi_n)}$ is a 
finite $($positive$)$ measure.
\end{prop}

One immediate implication of Proposition~\ref{ploc2.1} is the Leibniz rule, or derivation property,
\[
\ce_\chi(\varphi_1\varphi_2,\psi_1\psi_2)=
\ce_{\varphi_1\psi_1\chi}(\varphi_2,\psi_2)
+\ce_{\varphi_1\psi_2\chi}(\varphi_2,\psi_1)
+\ce_{\varphi_2\psi_1\chi}(\varphi_1,\psi_2)
+\ce_{\varphi_2\psi_2\chi}(\varphi_1,\psi_1)
\]
for the bilinear form  $\ce_\chi(\varphi,\psi)$ 
related to the truncated form $\ce_\chi$ by polarization and  a similar identity for the bilinear form $\ce(\varphi, \psi)$
associated with $\ce$.
(The latter identity is  formally obtained from the former by setting $\chi=\one_X$
and $\ce_{\one_X}=\ce$.)

\smallskip

Our  next application of the proposition is  a key identity related to the
estimate given by Condition~\ref{tloc1.1-2} of Theorem~\ref{tloc1.1}.

\begin{lemma} \label{lloc1} 
If $\varphi,\chi \in B(\ce)$ then
\begin{equation}
\ce_{\varphi^2}(\chi)
=\ce(\chi\,\varphi)-\ce(\varphi, \chi^2\varphi)
\;.
\label{eloc2.2}
\end{equation}
\end{lemma}
\proof\
 Let $\sigma_{ij}$ denote the representing measure of Proposition~\ref{ploc2.1} corresponding to the pair $\varphi, \chi$.
 Then
\begin{equation}
 \ce(\chi\,\varphi)=\int d\sigma_{11}(x_1,x_2)\,x_2^2+2\int d\sigma_{12}(x_1,x_2)\,x_1x_2
 +\int d\sigma_{22}(x_1,x_2)\,x_1^2
 \label{eloc2.20}
 \end{equation}
and
  \[
 \ce(\varphi, \chi^2\varphi)=\int d\sigma_{11}(x_1,x_2)\,x_2^2+2\int d\sigma_{12}(x_1,x_2)\,x_1x_2
 \;.
 \]
Therefore 
 \begin{equation}
\ce(\chi\,\varphi)- \ce(\varphi, \chi^2\varphi)=\int d\sigma_{22}(x_1,x_2)\,x_1^2
 \;.\label{eloc2.3}
 \end{equation}
Similarly, one calculates that 
 \begin{eqnarray}
\ce_{\varphi^2}(\chi)= \ce(\chi, \varphi^2\chi)- 2^{-1}\ce(\varphi^2,\chi^2)&=&\int d\sigma_{22}(x_1,x_2)\,x_1^2\
  \label{eloc2.4}
 \;
 \end{eqnarray}
Then (\ref{eloc2.2}) follows directly from (\ref{eloc2.3}) and (\ref{eloc2.4}).
\hfill$\Box$

\bigskip
The relevance of the identity (\ref{eloc2.2}) is that it formally identifies the energy measure corresponding to $\ce$ with a double commutator.
To illustrate this assume $\ce$ has a {\it carr\'e du champ} $\Gamma$.
Then it follows that 
\[
\ce_{\varphi^2}(\psi)=\int_X d\mu\,\varphi^2\,\Gamma(\psi)=(\varphi, \Gamma(\psi)\varphi)
\]
for $\varphi,\psi\in B(\ce)$.
Therefore  (\ref{eloc2.2}) gives the identification
\[
-2\,(\varphi, \Gamma(\psi)\varphi)=\ce(\varphi, \psi^2\varphi)-2\,\ce(\psi\varphi,\psi\varphi)+\ce(\psi^2\varphi, \varphi)
\]
for all  $\varphi,\psi\in B(\ce)$.
But if  $\psi\,\varphi\in D(H)$ for each $\varphi\in D(H)$ then 
\[
-2\,(\varphi, \Gamma(\psi)\varphi)=(H\varphi, \psi^2\varphi)-2\,(\psi\,\varphi,H\psi\,\varphi)+(\psi^2\varphi,H \varphi)
\]
which  is equivalent to the identification
\[
\Gamma(\psi)=-2^{-1}[\,[\,H,\psi\,],\,\psi\,]
\]
analogous to the situation for the Laplacian discussed in the introduction.

\smallskip

Next we need the following estimate.
\begin{lemma}\label{lloc2.10}
If $\varphi, \chi\in B(\ce)$ with $\chi\geq0$ then
\[
\ce_{\varphi^2}(\chi(1+\beta\chi)^{-1})\leq \ce_{(1+\beta\chi)^{-2}\varphi^2}(\chi)
\]
for all $\beta\geq0$.
\end{lemma}
\proof\
Again let $\sigma_{ij}$ denote the representing measure corresponding to the pair $\varphi, \chi$.
Then one calculates that 
\begin{eqnarray*}
\ce_{\varphi^2}(\chi(1+\beta\chi)^{-1})&=&\ce(\chi(1+\beta\chi)^{-1},\varphi^2\,\chi(1+\beta\chi)^{-1})-2^{-1}\ce(\chi^2(1+\beta\chi)^{-2},\varphi^2)\\[5pt]
&=&\int d\sigma_{22}(x_1,x_2)\,x_1^2\,(1+\beta x_2)^{-4}
\end{eqnarray*}
since the terms corresponding to $\sigma_{12}$ cancel.
Similarly
\begin{eqnarray*}
 \ce_{(1+\beta\chi)^{-2}\varphi^2}(\chi)&=&\ce(\chi,(1+\beta\chi)^{-2}\varphi^2\chi)-2^{-1}\ce(\chi^2,(1+\beta\chi)^{-2}\varphi^2)\\[5pt]
 &=&\int d\sigma_{22}(x_1,x_2)\,x_1^2\,(1+\beta x_2)^{-2}
 \end{eqnarray*}
 because the cross terms again cancel.
 (Since $\chi\geq0$ the $x_2$-integration is over the positive half axis.)
 Therefore the statement of the lemma follows immediately.
 \hfill$\Box$
 
 \bigskip

A locality estimate also gives  the following bounds.
\begin{lemma}\label{lloc2.2}
If $\varphi, \chi\in B(\ce)$  then 
\[
\ce(\chi\varphi)\leq(1+\delta)\,\ce_{\chi^2}(\varphi)+(1+\delta^{-1})\,\ce_{\varphi^2}(\chi)
\]
for all $\delta>0$.
If, in addition, $\psi\in B(\ce)$ then
\[
\ce_{\psi^2}(\chi\varphi)\leq(1+\delta)\,\ce_{\psi^2\chi^2}(\varphi)+(1+\delta^{-1})\,\ce_{\psi^2\varphi^2}(\chi)
\]
for all $\delta>0$
\end{lemma}
\proof\
It follows from (\ref{eloc2.20}) and the Cauchy--Schwarz inequality for the measure $\sigma_{ij}$   that
\begin{eqnarray*}
\ce(\chi\,\varphi)
 \leq(1+\delta^{-1})\int d\sigma_{11}(x_1,x_2)\,x_2^2+(1+\delta)\int d\sigma_{22}(x_1,x_2)\,x_1^2
  \;.
\end{eqnarray*}
But the second integral  on the right is equal to $\ce_{\chi^2}(\varphi)$ by (\ref{eloc2.4}).
Moreover, by interchanging $\chi$ and $\varphi$ one can identify the first integral with 
$\ce_{\varphi^2}(\chi)$.
The first statement of the lemma follows by combination of these observations.
The second statement follows by a similar calculation.
\hfill$\Box$

\bigskip
Next one has the following identity.
\begin{lemma}\label{lloc2.3}
If $\varphi, \chi\in B(\ce)$  then 
\[
\ce_{\varphi^2}(\chi^2)=4\,\ce_{\chi^2\varphi^2}(\chi)
\;.
\]
\end{lemma}

The proof follows by  a similar calculation
 to the derivation of  (\ref{eloc2.4}).
Alternatively it  follows directly from the Leibniz rule.

Finally we consider approximation of  functions in the domain $D(H)$ of the submarkovian operator $H$.
One has $D(H)\subseteq D(\ce)$ but  it is convenient to establish an explicit approximation,
 in the $D(\ce)$-graph norm,  of functions in $D(H)$  by functions in~$B(\ce)$.

\begin{lemma}\label{lloc3}
If $\varphi\in D(H)$ and $\varepsilon>0$ set $\varphi_\varepsilon=\varphi\,(1+\varepsilon\,\varphi^2)^{-1/2}$.
Then $\varphi_\varepsilon\in B(\ce)$, $\|\varphi_\varepsilon\|_2\leq\|\varphi\|_2$, $\ce(\varphi_\varepsilon)\leq
\ce(\varphi)$  and $\|\varphi-\varphi_\varepsilon\|_{D(\ce)}\to0$ as $\varepsilon\to0$.
Moreover, $\ce_{\xi}(\varphi_\varepsilon)\leq \ce_\xi(\varphi)$ for all $\xi\in B(\ce)_+$
and $\ce_{\xi}(\varphi_\varepsilon)\to \ce_\xi(\varphi)$ as $\varepsilon\to0$.
\end{lemma}
\proof\
First the boundedness property $\|\varphi_\varepsilon\|_2\leq\|\varphi\|_2$ is evident.
Secondly
\[
\|\varphi-\varphi_\varepsilon\|_2^2=\|\varphi\,(1-(1+\varepsilon\,\varphi^2)^{-1/2})\|_2^2
=\int_\Omega d\mu\,\varphi^2(1-(1+\varepsilon\,\varphi^2)^{-1/2})^2\leq \|\varphi\|_2^2
\;.
\]
Therefore   $\|\varphi-\varphi_\varepsilon\|_2\to0$ as $\varepsilon\to0$ by the Lebesgue dominated
 convergence theorem.

Thirdly,  the map $\varphi\mapsto \varphi_\varepsilon$ is a normal contraction
so $\ce(\varphi_\varepsilon)\leq\ce(\varphi)$ and $\ce_{\xi}(\varphi_\varepsilon)\leq \ce_\xi(\varphi)$
by the Markovian property of the form $\ce$ and the associated truncated functions.
The remaining convergence statements follow from the Anderssen representation, e.g.\ if
$\sigma^\varphi$ denotes the positive measure corresponding to $\varphi\in D(\ce)$ then
\[
\ce(\varphi-\varphi_\varepsilon)=\int d\sigma^\varphi(x)\,(1-(1+\varepsilon\,x^2)^{-3/2})^2\leq \ce(\varphi)
\;.
\]
Therefore $\ce(\varphi-\varphi_\varepsilon)\to0$ as $\varepsilon\to0$ by another application of dominated convergence.
Consequently,  $\varphi_\varepsilon\to \varphi$ as $\varepsilon\to0$ with respect to 
the $D(\ce)$-graph norm.
\hfill$\Box$

\bigskip

At this point we are prepared to 
 prove Theorem~\ref{tloc1.1}.
 The proof of  is in two steps.
First we establish the conclusion for  functions  $\varphi\in D(H)\cap B_c(\ce)$.
Secondly, we extend the result to all $\varphi\in D(H)$ by approximation.
The first step is straightforward but the second step is more complicated.
It involves simultaneous approximation of the Hardy function $\eta$ by bounded functions
of compact support and $\varphi$ by bounded functions in the form domain.

\medskip
\noindent{\bf Proof of Theorem~\ref{tloc1.1}}$\;$
Let  $\eta\in B_{\rm loc}(\ce)$ be  such that the $\eta$-Hardy inequality $\ce(\varphi)\geq (\eta\varphi, \eta\varphi)$ 
is satisfied for all $\varphi\in B_c(\ce)$.
Since $B(\ce)$ is an algebra it follows that $B_{\rm loc}(\ce)B_c(\ce)\subseteq B_c(\ce)$.
Therefore  $\ce(\eta\,\varphi)\geq (\eta^2\varphi, \eta^2\varphi)$ for all  $\varphi\in B_c(\ce)$.
Then it follows from Lemma~\ref{lloc1}  that 
\[
 \ce(\varphi, \eta^2\varphi)\geq (\eta^2\varphi, \eta^2\varphi)-\ce_{\varphi^2}(\eta)
\]
for all $\varphi\in B_c(\ce)$.
But   $\ce_{\varphi^2}(\eta)\leq \gamma\,(\eta^2\varphi,\eta^2\varphi)$ with $\gamma\in\langle0,1]$ 
by Condition~\ref{tloc1.1-2} of the theorem.
Therefore
\[
\ce(\varphi, \eta^2\varphi)\geq (1-\gamma)(\eta^2\varphi, \eta^2\varphi)
\]
for all $\varphi\in B_c(\ce)$.
If, however,  $\varphi\in D(H)\cap B_c(\ce)$ then $|\ce(\varphi, \eta^2\varphi)|\leq \|H\varphi\|_2\,\| \eta^2\varphi\|_2$ and one deduces that 
\[
\|H\varphi\|_2\,\| \eta^2\varphi\|_2\geq (1-\gamma)(\eta^2\varphi, \eta^2\varphi)=(1-\gamma)\|\eta^2\varphi\|_2^2
\;.
\]
Since $\gamma\in\langle0,1\rangle$ one
can divide by $\| \eta^2\varphi\|_2$ and square to obtain the Rellich inequality
\[
(H\varphi,H\varphi)=\|H\varphi\|_2^2\geq (1-\gamma)^2(\eta^2\varphi, \eta^2\varphi)
\]
for all $\varphi\in D(H)\cap B_c(\ce)$.

\smallskip

The derivation of the inequality for general $\varphi\in D(H)$ follows a similar line of reasoning
but is essentially more complicated.
As a preliminary note that  the  definition of the $\eta$-Hardy inequality  includes the condition $\eta\,D(\ce)\subseteq L_2(X)$
but the Rellich inequality requires that $\eta^2\,D(H)\subseteq L_2(X)$.
Since  $\eta\in  B_{\rm loc}(\ce)$ and $B(\ce)$ is an algebra one always has the weaker properties $\eta\, B_c(\ce)\subseteq L_2(X)$ and   $\eta^2\, B_c(\ce) \subseteq L_2(X)$. 
These relations do not depend on the validity of the Hardy or Rellich inequalities.
But the conditions $\eta\,D(\ce)\subseteq L_2(X)$ and $\eta^2\,D(H)\subseteq L_2(X)$
are much more stringent.
This is the  reason that the proof for general $\varphi\in D(H)$  is  complicated by  additional approximation arguments.

\smallskip

Fix $\varphi\in D(H)$ and set
 $\varphi_\varepsilon=\varphi\,(1+\varepsilon\,\varphi^2)^{-1/2}$ 
with $\varepsilon>0$.
Then
$\varphi_\varepsilon\in B(\ce)$, by Lemma~\ref{lloc3}.
Further set  $\eta_{\alpha,\beta}=\rho_\alpha\,\eta_\beta$ where the $\rho_\alpha\in B_c(\ce)$ satisfy 
Condition~\ref{tloc1.1-3} of the theorem and $\eta_\beta=\eta\,(1+\beta\eta)^{-1}$ with $\beta>0$. 
Then $\eta_\beta\in B(\ce)$, $\eta_{\alpha,\beta}\in B_c(\ce)$ and $0\leq\eta_{\alpha,\beta}\leq \eta_\beta\leq \eta$.
Moreover,
\begin{eqnarray}
(H\varphi, \eta_{\alpha,\beta}^{\;2}\,\varphi_\varepsilon)&=&\ce(\varphi, \eta_{\alpha,\beta}^{\;2}\,\varphi_\varepsilon)\nonumber\\[5pt]
&=&\ce(\varphi_\varepsilon, \eta_{\alpha,\beta}^{\;2}\,\varphi_\varepsilon)+\ce(\varphi-\varphi_\varepsilon, \eta_{\alpha,\beta}^{\;2}\,\varphi_\varepsilon)\nonumber\\[5pt]
&=&\ce(\eta_{\alpha,\beta}\,\varphi_\varepsilon)-\ce_{\varphi_\varepsilon^2}(\eta_{\alpha,\beta})
+\ce(\varphi-\varphi_\varepsilon, \eta_{\alpha,\beta}^{\;2}\,\varphi_\varepsilon)
\label{eloc201}
\end{eqnarray}
where the third step follows from Lemma~\ref{lloc1}.
In addition $\ce(\eta_{\alpha,\beta}\,\varphi_\varepsilon)\geq (\eta\,\eta_{\alpha,\beta}\,\varphi_\varepsilon,\eta\,\eta_{\alpha,\beta}\,\varphi_\varepsilon)$ by the $\eta$-Hardy inequality,
 Condition~\ref{tloc1.1-1}  
 of the theorem ,applied to $\eta_{\alpha,\beta}\,\varphi_\varepsilon$. 
Therefore
\begin{eqnarray}
(H\varphi, \eta_{\alpha,\beta}^{\;2}\,\varphi_\varepsilon)&\geq &
(\eta\,\eta_{\alpha,\beta}\,\varphi_\varepsilon,\eta\,\eta_{\alpha,\beta}\,\varphi_\varepsilon)-\ce_{\varphi_\varepsilon^2}(\eta_{\alpha,\beta})
+\ce(\varphi-\varphi_\varepsilon, \eta_{\alpha,\beta}^{\;2}\,\varphi_\varepsilon)\nonumber\\[5pt]
&\geq &
(\eta_{\alpha,\beta}\,\eta\,\varphi_\varepsilon,\eta_{\alpha,\beta}\,\eta\,\varphi_\varepsilon)
-\ce_{\varphi_\varepsilon^2}(\eta_{\alpha,\beta})
-\ce(\varphi-\varphi_\varepsilon)^{1/2}\, \ce(\eta_{\alpha,\beta}^{\;2}\,\varphi_\varepsilon)^{1/2}
\label{eloc21}
\end{eqnarray}
for all $\varphi\in D(H)$,  all $\beta,\varepsilon>0$ and all $\alpha$.
But
\[
|(H\varphi, \eta_{\alpha,\beta}^{\;2}\,\varphi_\varepsilon)|\leq \|H\varphi\|_2\,\|\eta_\beta^{\;2}\varphi_\varepsilon\|_2
\leq \|H\varphi\|_2\,\|\eta_\beta\,\eta\,\varphi_\varepsilon\|_2
\]
since $\eta\,\varphi_\varepsilon\in L_2(X)$ by the $\eta$-Hardy inequality.
Combining this estimate with (\ref{eloc21}) and taking a limit over $\alpha$ then gives
\begin{eqnarray}
\|H\varphi\|_2\,\|\eta_\beta\,\eta\,\varphi_\varepsilon\|_2\geq \|\eta_\beta\,\eta\,\varphi_\varepsilon\|_2^2
&-&\limsup_\alpha\ce_{\varphi_\varepsilon^2}(\eta_{\alpha,\beta})\nonumber\\[0pt]
&-&\ce(\varphi-\varphi_\varepsilon)^{1/2}\,\limsup_\alpha\ce(\eta_{\alpha,\beta}^{\;2}\,\varphi_\varepsilon)^{1/2}
\label{eloc210}
\end{eqnarray}
for all $\beta,\varepsilon>0$.
Here we have used Condition~\ref{tloc1.1-3} to deduce that $\eta_{\alpha,\beta}$ converges on $L_2(X)$ to $\eta_\beta$.
It is important at this point that $\eta\,\varphi_\varepsilon\in L_2(X)$ by the $\eta$-Hardy inequality and $\eta_\beta\in L_\infty(X)$.
Therefore $\eta_\beta\,\eta\,\varphi_\varepsilon\in L_2(X)$.

Next consider the second term on the right hand side of (\ref{eloc210}).
Since $\eta_{\alpha,\beta}=\rho_\alpha\,\eta_\beta$  it follows from Lemma~\ref{lloc2.2}  that 
\[
\ce_{\varphi_\varepsilon^2}(\eta_{\alpha,\beta})
\leq (1+\delta)\,\ce_{\rho_\alpha^2\varphi_\varepsilon^2}(\eta_\beta)+(1+\delta^{-1})\,\ce_{\eta_\beta^2\varphi_\varepsilon^2}(\rho_\alpha)
\]
for all $\delta>0$.
Now we apply Lemma~\ref{lloc2.10}, with $\varphi$ replaced by $\rho_\alpha\varphi_\varepsilon$ and
$\chi$ replaced by $\eta$, and   Condition~\ref{tloc1.1-2} of the theorem,
with $\varphi$ replaced by $(1+\beta\eta)^{-1}\rho_\alpha\varphi_\varepsilon$, to the first term.
One finds
\begin{eqnarray*}
\ce_{\rho_\alpha^2\varphi_\varepsilon^2}(\eta_\beta)&\leq&\ce_{(1+\beta\eta)^{-2}\rho_\alpha^2\varphi_\varepsilon^2}(\eta)\\[5pt]
&\leq&\gamma\,(\eta^2(1+\beta\eta)^{-1}\rho_\alpha\varphi_\varepsilon,\eta^2(1+\beta\eta)^{-1}\rho_\alpha\varphi_\varepsilon)
\leq\gamma\,\|\eta_{\beta}\,\eta\,\varphi_\varepsilon\|_2^2
\end{eqnarray*}
where we again have $\eta\,\varphi_\varepsilon\in L_2(X)$ by the $\eta$-Hardy inequality and $\eta_\beta\in L_\infty(X)$.
Note that these steps are valid because $\rho_\alpha$ has compact support.
Therefore one now has
\begin{eqnarray*}
\ce_{\varphi_\varepsilon^2}(\eta_{\alpha,\beta})
\leq (1+\delta)\,\gamma\,\|\eta_{\beta}\,\eta\,\varphi_\varepsilon\|_2^2+(1+\delta^{-1})\,\ce_{\eta_\beta^2\varphi_\varepsilon^2}(\rho_\alpha)\end{eqnarray*}
for all $\delta>0$.
Hence 
\[
\limsup_\alpha\ce_{\varphi_\varepsilon^2}(\eta_{\alpha,\beta})\leq  (1+\delta)\,\gamma\,\|\eta_{\beta}\,\eta\,\varphi_\varepsilon\|_2^2
\]
for all $\delta>0$ since $\limsup_\alpha\ce_{\eta_\beta^2\varphi_\varepsilon^2}(\rho_\alpha)=0$
by   Condition~\ref{tloc1.1-3} of the theorem.
Combining this estimate with (\ref{eloc210}) and taking the limit of $\delta\to0$ one then obtains the bounds
\begin{eqnarray}
\|H\varphi\|_2\,\|\eta_\beta\,\eta\,\varphi_\varepsilon\|_2&\geq&(1-\gamma)\, \|\eta_\beta\,\eta\,\varphi_\varepsilon\|_2^2
-\ce(\varphi-\varphi_\varepsilon)^{1/2}\,\limsup_\alpha\ce(\eta_{\alpha,\beta}^{\;2}\,\varphi_\varepsilon)^{1/2}
\label{eloc211}
\end{eqnarray}
for all $\beta,\varepsilon>0$.
Our next aim it to prove that $\limsup_\alpha\ce( \eta_{\alpha,\beta}^{\;2}\,\varphi_\varepsilon)$ is bounded uniformly in~$\varepsilon$.

First consider $\ce(\eta_{\alpha,\beta}\,\varphi_\varepsilon)$.
It follows that 
\[
\ce(\eta_{\alpha,\beta}\,\varphi_\varepsilon)\leq 2\,\ce_{\eta_{\alpha,\beta}^{\;2}}(\varphi_\varepsilon)
+2\,\ce_{\varphi_\varepsilon^{\,2}}(\eta_{\alpha,\beta})
\leq 2\,\|\eta_\beta\|_\infty^2\,\ce(\varphi_\varepsilon)
+2\,\ce_{\varphi_\varepsilon^{\,2}}(\eta_{\alpha,\beta})
\]
where the first estimate uses  Lemma~\ref{lloc2.2} and the second uses $\ce_{\chi^2}(\psi)\leq \|\chi\|_\infty^2\,\ce(\psi)$
and $\|\eta_{\alpha,\beta}\|_\infty\leq \|\eta_\beta\|_\infty$.
Then another application of  Lemma~\ref{lloc2.2}  gives
\[
\ce_{\varphi_\varepsilon^{\,2}}(\eta_{\alpha,\beta})\leq 2\,
\ce_{\rho_\alpha^2\,\varphi_\varepsilon^{\,2}}(\eta_{\beta})
+2\,\ce_{\eta_\beta^{\;2}\,\varphi_\varepsilon^{\,2}}(\rho_{\alpha})
\;.
\]
Combining these bounds
 one obtains
\[
\ce( \eta_{\alpha,\beta}\,\varphi_\varepsilon)\leq 
2\,\|\eta_\beta\|_\infty^2\,\ce(\varphi_\varepsilon)+
4\,\ce_{\rho_\alpha^2\,\varphi_\varepsilon^{\,2}}(\eta_{\beta})
+4\,\ce_{\eta_\beta^{\;2}\,\varphi_\varepsilon^{\,2}}(\rho_{\alpha})
\;.\]
But $\ce(\varphi_\varepsilon)\leq \ce(\varphi)$ by Lemma~\ref{lloc3}.
Moreover,
\begin{equation}
\ce_{\rho_\alpha^2\,\varphi_\varepsilon^{\,2}}(\eta_{\beta})
\leq \ce_{(1+\beta\eta)^{-2}\rho_\alpha^2\,\varphi_\varepsilon^{\;2}}(\eta)
\leq \gamma\,(\eta_\beta\,\eta\,\rho_\alpha\,\varphi_\varepsilon,\eta_\beta\,\eta\,\rho_\alpha\,\varphi_\varepsilon)
\leq \gamma\,\|\eta_\beta\,\eta\,\varphi\|_2^2
\label{eloc212}
\end{equation}
by another application of Lemma~\ref{lloc2.10} and Condition~\ref{tloc1.1-2} of the theorem.
Again this is valid since $\rho_\alpha$ has compact support.
Therefore combination of these last two estimates gives
\[
\limsup_\alpha\ce( \eta_{\alpha,\beta}\,\varphi_\varepsilon)\leq 
2\,\|\eta_\beta\|_\infty^2\,\ce(\varphi)+4\,\gamma\,\|\eta_\beta\,\eta\,\varphi\|_2^2
\]
since $\limsup_\alpha \ce_{\eta_\beta^{\;2}\,\varphi_\varepsilon^{\,2}}(\rho_{\alpha})=0$
by Condition~\ref{tloc1.1-3}.
Note that this last bound is uniform in $\varepsilon$.

Next by Lemma~\ref{lloc2.2} one has
\begin{eqnarray*}
\ce( \eta_{\alpha,\beta}^{\;2}\,\varphi_\varepsilon)\leq2\,\ce_{\eta_{\alpha,\beta}^{\;2}}(\eta_{\alpha,\beta}\,\varphi_\varepsilon)
+2\,\ce_{\eta_{\alpha,\beta}^{\;2}\varphi_\varepsilon^{\;2}}(\eta_{\alpha,\beta})
\leq2\,\|\eta_{\beta}\|_\infty^2\,\ce(\eta_{\alpha,\beta}\,\varphi_\varepsilon)
+2\,\ce_{\eta_{\alpha,\beta}^{\;2}\varphi_\varepsilon^{\;2}}(\eta_{\alpha,\beta})
\end{eqnarray*}
and arguing as in the last paragraph
\begin{eqnarray*}
\ce_{\eta_{\alpha,\beta}^{\;2}\varphi_\varepsilon^{\;2}}(\eta_{\alpha,\beta})
\leq 2\,\ce_{\rho_\alpha^{\;2}\eta_{\alpha,\beta}^{\;2}\varphi_\varepsilon^{\;2}}(\eta_{\beta})
+2\,\ce_{\eta_\beta^{\;2}\eta_{\alpha,\beta}^{\;2}\varphi_\varepsilon^{\;2}}(\rho_{\alpha})
\leq 2\,\ce_{\rho_\alpha^{\;2}\eta_{\beta}^{\;2}\varphi_\varepsilon^{\;2}}(\eta_{\beta})
+2\,\ce_{\eta_\beta^{\;4}\varphi_\varepsilon^{\;2}}(\rho_{\alpha})
\;.
\end{eqnarray*}
Therefore combining these estimates one finds
\[
\limsup_\alpha\ce( \eta_{\alpha,\beta}^{\;2}\,\varphi_\varepsilon)
\leq 2\,\|\eta_{\beta}\|_\infty^2\,\limsup_\alpha\ce( \eta_{\alpha,\beta}\,\varphi_\varepsilon)
+4\,\limsup_\alpha \ce_{\rho_\alpha^{\;2}\eta_{\beta}^{\;2}\varphi_\varepsilon^{\;2}}(\eta_{\beta})
\;.
\]
The first term on the right is bounded uniformly in $\varepsilon$ by the previous argument  and
\[
\ce_{\rho_\alpha^{\;2}\eta_{\beta}^{\;2}\varphi_\varepsilon^{\;2}}(\eta_{\beta})
\leq \gamma\,(\eta_\beta^{\;2}\,\eta\,\rho_\alpha\,\varphi_\varepsilon, \eta_\beta^{\;2}\,\eta\,\rho_\alpha\varphi_\varepsilon)
\leq \gamma\,\|\eta_\beta^{\;2}\,\eta\,\varphi\|_2^2
\]
by the estimate (\ref{eloc212}) with $\varphi_\varepsilon$ replaced by $\eta_\beta\,\varphi_\varepsilon$.
Therefore $\limsup_\alpha\ce( \eta_{\alpha,\beta}^{\;2}\,\varphi_\varepsilon)$ is bounded uniformly in $\varepsilon$.

Now one can take the limit $\varepsilon\to0$ in  (\ref{eloc211}).
Since 
 $\|\varphi-\varphi_\varepsilon\|_{D(\ce)}\to0$ by Lemma~\ref{lloc3}, and $\eta\,\varphi\in L_2(X)$ by the $\eta$-Hardy inequality,  
 it follows that 
\begin{eqnarray*}
\|H\varphi\|_2\,\|\eta_\beta\,\eta\,\varphi\|_2
&=&\lim_{\varepsilon\to0}\|H\varphi\|_2\,\|\eta_\beta\,\eta\,\varphi_\varepsilon\|_2\\[5pt]
&\geq& (1-\gamma)\,\limsup_{\varepsilon\to0} \|\eta_\beta\,\eta\,\varphi_\varepsilon\|_2^2
=(1-\gamma)\, \|\eta_\beta\,\eta\,\varphi\|_2^2
\;.
\end{eqnarray*}
Therefore 
\[
\|H\varphi\|_2\,\|\eta_\beta\,\eta\,\varphi\|_2\geq  (1-\gamma)\,\|\eta_{\beta}\,\eta\,\varphi\|_2^2
\;.
\]
Since $\gamma<1$ by assumption one deduces that
 \[
 \|H\varphi\|_2^2\geq (1-\gamma)^2\,\|\eta_{\beta}\,\eta\,\varphi\|_2^2
 \;.
 \]
 Finally since the left hand side is independent of $\beta$ one concludes by dominated convergence that
 $\eta^2 \varphi\in D(H)$ and 
  \[
 \|H\varphi\|_2^2\geq (1-\gamma)^2\,\|\eta^2\varphi\|_2^2
\]
for all $\varphi\in D(H)$.\hfill$\Box$

 \bigskip
 
 Condition~\ref{tloc1.1-3} of Theorem~\ref{tloc1.1} has a different nature to the first two conditions since it is independent of $\eta$.
 It is related to the existence 
 of an approximate identity in the $D(\ce)$-graph norm.
 In particular if
 there is a net $\{\rho_\alpha\}$ with 
$\rho_\alpha\in B_c(\ce)$ such that 
$0\leq\rho_\alpha\leq 1$, 
$\rho_\alpha \,D(\ce)\subseteq D(\ce)$ for all $\alpha$ and 
\begin{equation}
\lim_\alpha\|(\rho_\alpha-\one_X)\varphi\|_{D(\ce)}=0
\label{ploc3.11}
\end{equation}
for all $\varphi\in D(\ce)$ then  Condition~\ref{tloc1.1-3}  is satisfied.
This follows because 
\begin{eqnarray*}
\ce_{\varphi^2}(\rho_\alpha)&=&\ce(\rho_\alpha\varphi)-\ce(\varphi, \rho_\alpha^2\varphi)\\[5pt]
&=&\ce((\rho_\alpha-\one_X)\varphi)-\ce(\varphi, (\rho_\alpha-\one_X)^2\varphi)
\end{eqnarray*}
where the first step follows from Lemma~\ref{lloc1}  and the second by direct calculation.
But $|\ce(\varphi, (\rho_\alpha-\one_X)^2\varphi)|^2\leq \ce(\varphi)\,\ce((\rho_\alpha-\one_X)^2\varphi)$
by the Cauchy--Schwarz inequality.
Then  it follows from (\ref{ploc3.11}) and the  uniform boundedness principle  that there is an $M>0$ such that
 $\ce((\rho_\alpha-\one_X)^2\varphi)\leq M\,\ce((\rho_\alpha-\one_X)\varphi)$ for all $\alpha$.
 Therefore $\lim_\alpha\ce_{\varphi^2}(\rho_\alpha)=0$.
 In fact there is a weak converse to this statement: if Condition~\ref{tloc1.1-3} of Theorem~\ref{tloc1.1}
 is satisfied for all $\varphi\in B(\ce)$ then (\ref{ploc3.11}) is valid for all $\varphi\in B(\ce)$.

\section{Operators on domains}\label{S4}

In this section we give the proof of Theorem~\ref{tloc1.2}.
To be more precise we deduce the theorem as a corollary of Theorem~\ref{tloc1.1}.
The principal difficulty is to verify Condition~\ref{tloc1.1-3} of the latter theorem,
the existence of a suitable approximate identity.
This is the critical property used in Section~\ref{S2} to extend the Rellich inequality
from functions with compact support 
to the full domain of the submarkovian operator.

\medskip

\noindent{\bf Proof of Theorem~\ref{tloc1.2}}$\;$
First the Hardy inequality (\ref{eloc1.6}) which is the principal assumption of the theorem 
can be reformulated as the $\eta$-Hardy inequality (\ref{eloc1.1}) by choosing
  $\eta$ as the positive square root of  $a_1\,c_\Omega\,d_\Omega^{\;-2}$.
Then $\eta\in W^{1,2}_{\rm loc}(\Omega)=B_{\rm loc}(h)$ and Condition~\ref{tloc1.1-1} of  
Theorem~\ref{tloc1.1} is verified.
Moreover, $h_{\varphi^2}(\eta)=(\varphi, \Gamma(\eta)\varphi)$ 
where the {\it carr\'e du champ} $\Gamma$ is given by $\Gamma(\eta)=c_\Omega\,|\nabla\eta|^2$.
But $4\,\eta^2\,\Gamma(\eta)=\Gamma(\eta^2)$.
Therefore a straightforward calculation, using $|\nabla d_\Omega|=1$, gives
\[
\Gamma(\eta)=
a_1\,c_\Omega^{\;2}\,d_\Omega^{\;-4}|1-d_\Omega \,c'_\Omega/(2\,c_\Omega)|^2
\leq (\nu/a_1)\,\eta^4
\]
where $c'_\Omega=c'\circ d_\Omega$. 
The last step follows since $\delta\wedge\delta'\leq (s\,c'(s)/c(s))\leq \delta\vee\delta'$
as a consequence of the identity $s\,c'(s)/c(s)=(\delta+\delta'\,s)/(1+s)$.
Hence one concludes that 
\[
h_{\varphi^2}(\eta)=(\varphi, \Gamma(\eta)\varphi)\leq\gamma\,(\eta^2\varphi,\eta^2\varphi)
\]
for all $\varphi\in B_c(h)$ with $\gamma=\nu/a_1$.
Then it follows from the first paragraph of the proof of Theorem~\ref{tloc1.1} that if $\gamma<1$ then 
\[
(H\varphi, H\varphi)\geq (1-\gamma)^2\,(\eta^2\varphi, \eta^2\varphi)
=a_1^2\,(1-\gamma)^2 \,(c_\Omega\, d_\Omega^{\;-2}\varphi,c_\Omega\, d_\Omega^{\;-2}\varphi)
\]
for all $\varphi\in B_c(h)$ and in particular for all $\varphi\in C_c^\infty(\Omega)$.
Thus the Rellich inequality (\ref{eloc1.7}) is satisfied for all $\varphi\in B_c(h)$ with $a_2=a_1^2\,(1-\gamma)^2=(a_1-\nu)^2$.
This is the elementary part of the proof.
The difficulty lies in extending the Rellich inequality to all $\varphi \in D(H)$.
This is achieved by the construction of an approximate identity satisfying Condition~\ref{tloc1.1-3} of Theorem~\ref{tloc1.1}.
The construction is an adaptation of a key idea of Agmon (see \cite{Agm1}, Chapter~1), the introduction of an alternative metric.

\medskip

The Euclidean metric on $\Omega$ is usually defined by a shortest path algorithm but it is also given by  the equivalent definition 
\[
d(x\,;y)=\sup\{|\psi(x)-\psi(y)|: \psi\in W_{\rm \loc}^{1,\infty}(\Omega),\, |\nabla\psi|\leq 1\,\}
\;.
\]
The  Euclidean distance from $x\in\Omega$  to the  measurable subset $A\subseteq \Omega$ is  then defined   by 
$d(x\,;A)=\inf_{y\in A}d(x\,;y)$ and  $d_\Omega$ 
is given by
 \[
 d_\Omega(x)=\sup\{d(x\,; \Omega\backslash K):K \mbox{  is a compact subset of } \Omega\}
 \;.
 \]
Now we use these definitions as the model for introducing an alternative metric.

Define the distance $d_2(\,\cdot\,;\,\cdot\,)$  by
\[
d_2(x\,;y)=\sup\{|\psi(x)-\psi(y)|: \psi\in W^{1,\infty}_{\rm loc}(\Omega),\, d_\Omega^{\;2}\,|\nabla\psi|^2\leq 1\,\}
\;.
\]
Then the  $d_2$-distance to the    measurable subset $A$  is given  by 
$d_2(x\,;A)=\inf_{y\in A}d(x\,;y)$ and the corresponding distance to the boundary by
 \[
 d_{2;\Omega}(x)=\sup\{d_2(x\,; \Omega\backslash K):K \mbox{  is a compact subset of } \Omega\}
 \;.
 \]
The functions $d_\Omega$ and  $d_{2;\Omega}$ are  Lipschitz and  satisfy
 $ |\nabla d_{\Omega}|^2=1$ and  $d_\Omega^{\;2}\, |\nabla d_{2;\Omega}|^2=1$ almost everywhere.
 
 The motivation for the introduction of $d_2(\,\cdot\,;\,\cdot\,)$ is the following.
 
 \begin{lemma}\label{lloc4.1}
 The metric space $(\Omega, d_2(\,\cdot\,;\,\cdot\,))$ is complete.
 \end{lemma}
 \proof\
 It follows from  Lemma~A1.2  in Appendix~A of Agmon's lecture notes \cite{Agm1} 
 that the space is complete if and only if $ d_{2;\Omega}(x)=\infty$ for one $x\in\Omega$ or,  equivalently,
  for all $x\in\Omega$.
  The latter equivalence is  a simple application of the triangle inequality.
  Therefore it suffices to  argue that the $d_2$-distance to the boundary is infinite.
  
  Let $\Omega_\delta=\{x\in\Omega: d_\Omega(x)<\delta\}$ for all $\delta>0$.
Fix $0<\delta_1<\delta_2$.
Then  introduce  $ \psi\in W^{1,\infty}_{\rm loc}(\Omega)$ such that 
\[
\psi(x)
= \left\{ \begin{array}{ll}
   -\log (\delta_1/\delta_2) \;\;\;\;& \mbox{if } x\in\Omega_{\delta_1}\; ,  \\[5pt]
    -\log (d_\Omega(x)/\delta_2)\;\;\;\;& \mbox{if } x \in\overline\Omega_{\delta_2}\backslash\overline\Omega_{\delta_1} \; , \\[5pt]
   \hspace{1cm}0 \;\;\;\;& \mbox{if } x\in \Omega\backslash\overline\Omega_{\delta_2} \;.          \end{array} \right.
\]
Since  $ d_\Omega^{\;2}\,|\nabla (\log d_\Omega)|^2=|\nabla d_\Omega|^2\leq1$
it follows that 
\[
d_2(x\,;y)\geq |\psi(x)-\psi(y)|=|\log (d_\Omega(x)/\delta_2)-\log (\delta_1/\delta_2)|=\log (d_\Omega(x)/\delta_1)
\]
for all $ x \in\overline\Omega_{\delta_2}\backslash\overline\Omega_{\delta_1}$ and $y\in\Omega_{\delta_1}$.
Therefore in the limit $\delta_1\to0$ one deduces that $d_{2;\Omega}(x)=\infty$ for all $x\in\Omega$.
\hfill$\Box$

\bigskip

The conclusion of the foregoing argument can be rephrased as follows.
\begin{cor}\label{cloc4.1}
 For each $m>0$ one has $\Omega=\{x\in\Omega: d_{2;\Omega}(x)>m\}$.
\end{cor}

The completeness property of Lemma~\ref{lloc4.1}  is crucial for the verification of Condition~\ref{tloc1.1-3} of Theorem~\ref{tloc1.1}.
The second crucial feature is the Hardy inequality.
But the following reasoning is not restricted to the forms and operators covered by Theorem~\ref{tloc1.2}.
The only aspect of the Hardy inequality of relevance is the property $c_\Omega^{\;1/2}d_\Omega^{\;-1}D(h)\subseteq L_2(\Omega)$.
We now construct a sequence of  $\rho_n$  satisfying the appropriate properties of an approximate identity  by the reasoning of Agmon \cite{Agm1} in   the proof of his Theorem~1.5. 

\begin{prop}\label{ploc4.1}
 Assume that $h$ satisfies the Hardy inequality $(\ref{eloc1.6})$.
 Then there exists a sequence $\rho_n\in D(h)$ with compact support such that $0\leq \rho_n\leq 1$ and 
 \[
\lim_{n\to\infty}\|(\rho_n-\one_\Omega)\varphi\|_2=0 \;\;\;\;\;{\rm and}\;\;\;\;\; \lim_{n\to\infty}h_{\varphi^2}(\rho_n)=0
 \] 
 for all $\varphi\in D(h)$.
 \end{prop}
 \proof\
 First  let $\{\Omega_n\}_{n\geq1}$ be a compact exhaustion sequence of $\Omega$, i.e.\
the $\Omega_n$ are relatively compact  open subsets of $\Omega$ with $\overline\Omega_n\subset \Omega_{n+1}$
such that  $\Omega=\bigcup_{n\geq1}\Omega_n$.
Secondly, fix $m>0$ and define $\rho$ by $\rho(t)=(t/m)\wedge 1$.
Further define $\rho_n$ by
$\rho_n(x)=\rho(d_2(x\,;\Omega\backslash \Omega_n))$.
Then
\[
\rho_n(x)
= \left\{ \begin{array}{ll}
  \hspace{1cm} 0 \;\;\;\;& \mbox{if } x\in\Omega\backslash \Omega_n\; ,  \\[5pt]
   d_2(x\,;\Omega\backslash \Omega_n)/m \;\;\;\;& \mbox{if } x \in \Omega_n \mbox{  and  } d_2(x\,;\Omega\backslash \Omega_n)\leq m \;,
   \\[5pt]          
   \hspace{1cm}1& \mbox{if } x \in \Omega_n \mbox{  and  } d_2(x\,;\Omega\backslash \Omega_n)>m\;.
\end{array} \right.
\]
It follows immediately from this definition that the  $\rho_n$ have compact support.
In particular $\supp\rho_n\subseteq \overline\Omega_n$.
Moreover,  $0\leq \rho_n\leq1$.
In addition the pointwise  limit of the $\rho_n$ as $n\to\infty$ is equal 
to the identity on the set of $x$ for which $d_{2;\Omega}(x)>m$.
But this set is equal to $\Omega$  by  Corollary~\ref{cloc4.1}.
Thus  the $\rho_n$ converge pointwise to $\one_\Omega$ as $n\to\infty$.
Hence they also converge to the identity  strongly on $L_2(\Omega)$.
Finally the $\rho_n$ are Lipschitz continuous.
This follows by first noting that 
\[
|\rho(s)-\rho(t)|\leq m^{-1}|s-t|
\;.
\]
Therefore
\begin{eqnarray*}
|\rho_n(x)-\rho_n(y)|&\leq& m^{-1}| d_2(x\,;\Omega\backslash \Omega_n)- d_2(y\,;\Omega\backslash \Omega_n)|\\[5pt]
&\leq& m^{-1}d_2(x\,;y)
\leq C\,m^{-1}|x-y|
\;.
\end{eqnarray*}
The second inequality follows from the triangle inequality and is valid for all $x,y\in\Omega$.
The third inequality follows from the definition of $d_2(\,\cdot\,;\,\cdot\,)$ and is valid locally, i.e.\ it is valid  in a neighbourhood of 
 each  fixed point $x_0\in\Omega$ with  the value of $C$ depending  on $x_0$.
Consequently, since the $\rho_n$ are Lipschitz functions and  $|\rho_n(x)-\rho_n(y)|\leq  m^{-1}d_2(x\,;y)$
for all $x, y\in \Omega$ it  follows from  the eikonal inequality (see \cite{Agm1}, Theorem~1.4 (ii))  that 
\[
d_\Omega^{\;2}\,|\nabla\!\rho_{n}|^2\leq  m^{-2}
\;.
\]
This is the critical inequality since it gives
\[
\Gamma(\rho_n)=c_\Omega\,|\nabla\!\rho_{n}|^2\leq m^{-2} \, c_\Omega\,d_\Omega^{\;-2}
\;.
\]
This estimate corresponds to (1.23) in Agmon's notes.
Next note that  the derivatives of $\rho_n$ have support in the set $\Omega_{m,n}= \{x \in \Omega_n: d_2(x\,;\Omega\backslash \Omega_n)\leq m\}$
and  $|\Omega_{m,n}|\to0$ as $n\to\infty$ by Corollary~\ref{cloc4.1}. 
Moreover, if $K$ is a compact subset of $\Omega$ then $K\subset \Omega_n$ for all sufficiently large $n$.
Up to this point we have not used the Hardy inequality (\ref{eloc1.6}).
But it follows from this inequality that  if $\varphi\in D(h)$ then $\varphi\in D(c_\Omega^{\;1/2} d_\Omega^{\;-1})$.
Therefore $\psi=c_\Omega^{\;1/2}d_\Omega^{\;-1}\varphi\in L_2(\Omega)$.
Hence 
\[
\int_\Omega\,\Gamma(\rho_n)\,|\varphi|^2\leq m^{-2}\int_{\Omega_{m,n}}\!\!\!c_\Omega \,d_\Omega^{\;-2}\,|\varphi|^2=
 m^{-2}\int_{\Omega_{m,n}}\!\! |\psi|^2
 \;.
 \]
 Since $\psi\in L_2(\Omega)$ it follows directly from this estimate that
  \[
 \lim_{n\to\infty}\int_\Omega\,\Gamma(\rho_n)\,|\varphi|^2=0
 \;.
 \]
 Thus, in the earlier notation,   $h_{\varphi^2}(\rho_n)\to0$ as $n\to\infty$.
  \hfill$\Box$

\bigskip

\noindent{\bf Proof of Theorem~\ref{tloc1.2} continued}
The proof of the theorem is now a corollary of Theorem~\ref{tloc1.1}.
Condition~\ref{tloc1.1-1} of the theorem is valid by assumption of the Hardy inequality (\ref{eloc1.6}), 
Condition~\ref{tloc1.1-2} was verified at the beginning of the proof with $\gamma=\nu/a_1$ and 
Condition~\ref{tloc1.1-3} now follows from Proposition~\ref{ploc4.1}.
Therefore the Rellich inequality (\ref{eloc1.7}) follows for all $\varphi\in D(H)$
with $a_2=(a_1-\nu)^2$
whenever $\nu<a_1$.\hfill$\Box$

\section{Applications and illustrations}\label{S5}

In this section we give two illustrations of Theorems~\ref{tloc1.1} and \ref{tloc1.2}.
First we give a direct  application of  the latter theorem with $\Omega=\Ri^d\backslash\{0\}$.
The application  requires establishing the validity  of the Hardy inequality (\ref{eloc1.6}), 
calculating the Hardy constant $a_1$ and  verifying the condition $a<\nu_1$.

\begin{exam}\label{exloc4.1}$\;$ Let  $\Omega=\Ri^d\backslash\{0\}$.
Then $\partial\Omega=\{0\}$ and $d_\Omega(x)=|x|$.
The Dirichlet form $h$, submarkovian operator $H$, cofficient function $c$ etc. are defined in Theorem~\ref{tloc1.2}.
In particular $c(s)=s^\delta(1+s)^{\delta'-\delta}$ with $\delta,\delta'\geq0$ and $c_\Omega(x)=c(|x|)$.

\begin{obs}\label{oloc5.1}
If $d+(\delta\wedge\delta')-2>0$ then the Hardy inequality
\begin{equation}
h(\varphi)\geq a_1 \,(c_\Omega^{\,1/2}d_\Omega^{\;-1}\varphi, c_\Omega^{\,1/2}d_\Omega^{\;-1}\varphi)
\label{eloc5.0}
\end{equation}
is valid for all $\varphi\in D(h)$ with $a_1=(d+(\delta\wedge\delta')-2)^2/4$.
The value of $a_1$ is optimal.
\end{obs}
\proof\
First note that 
\begin{eqnarray}
h(\varphi)-\lambda\,(\varphi, (\divv c_\Omega\,\chi)\varphi)&+&\lambda^2\,(\varphi,  c_\Omega\,\chi^2\,\varphi)
\nonumber \\[5pt]
&&\hspace{5mm}=\left((\nabla+\lambda \,\chi)\varphi,c_\Omega \,(\nabla+\lambda\,\chi)\varphi\right)\geq0
\label{eloc5.1}
\end{eqnarray}
for all $\varphi\in C_c^\infty(\Omega)$ and $\chi=(\chi_1,\ldots,\chi_d)$ with $\chi_k\in W^{1,\infty}_{\rm loc}(\Omega)$.
Now  choose $\chi=(\nabla d_\Omega)\, d_\Omega^{\;-1}$.
Thus $\chi(x)=x\,|x|^{-2}$.
Since $|\nabla d_\Omega|=1$ it follows that $  c_\Omega\,\chi^2=c_\Omega\,d_\Omega^{\;-2}$.
Moreover,
\[
\divv(c_\Omega\,\chi)=(d-2+c'_\Omega\,d_\Omega/c_\Omega)\,c_\Omega\,d_\Omega^{\;-2}
\geq (d+(\delta\wedge\delta')-2)\,c_\Omega\,d_\Omega^{\;-2}
\]
since $c'(s)\,s/c(s)\geq \delta\wedge\delta'$.
Then one deduces from (\ref{eloc5.1}) that 
\begin{eqnarray}
h(\varphi)-2\,\lambda\,b\,(\varphi, c_\Omega\,d_\Omega^{\,-2}\varphi)+\lambda^2\,(\varphi, c_\Omega\,d_\Omega^{\,-2}\varphi)\geq0
\label{eloc5.100}
\end{eqnarray}
for all $\lambda >0$ where $b=(d+(\delta\wedge\delta')-2)/2$.
Then if $b>0$ one can choose $\lambda=b$  and   the Hardy inequality 
follows, with $a_1=b^2=(d+(\delta\wedge\delta')-2)^2/4$,
for all $\varphi \in C_c^\infty(\Omega)$ and  then by continuity for all $\varphi\in D(h)$.

The optimality of $a_1$ follows by variation a standard argument (see, for example, \cite{BEL} Chapter~1).
First  let $a$ denote the optimal value of the constant for a Hardy inequality of the form (\ref {eloc5.0}).
Secondly, set $a(\delta) = ((d+\delta-2)/2)^2$.
Then it follows from (\ref{eloc5.0}) that $a\geq a(\delta\wedge \delta')$.  
Therefore it suffices to prove the identical upper bound. 
But then it is sufficient  to prove that $a(\delta)$ and $a(\delta')$ are both upper bounds since 
 $a(\delta\wedge \delta')=a(\delta)\wedge a(\delta')$.
The first bound is established by an estimate at the origin and the second by a similar estimate at infinity.

The estimate at the origin is obtained by examining  functions $\varphi_\alpha=d_\Omega^{\;-\alpha}\xi$, $\alpha>0$, where $\xi$ has support in  a small neighbourhood of the origin.
Then $c_\Omega\,d_\Omega^{\;-2}\,|\varphi_\alpha|^2$ is integrable if $\alpha <(d+\delta-2)/2$.
Choosing $\alpha=(d+\delta-2-\varepsilon)/2$ with $\varepsilon>0$ small
one can arrange that $h(\varphi_\alpha)/\|c_\Omega^{\,1/2}d_\Omega^{\,-1}\varphi_\alpha\|_2^2\leq \alpha^2$
and  in the limit $\varepsilon\to0$ one concludes that $a\leq a(\delta)$.
Here the property $\lim_{s\to0}c(s)\,s^{-\delta}=1$ is important.
The estimate at infinity is similar. 
One now chooses $\varphi_\alpha$ with support in the complement of  a large ball centred at the origin and equal to $d_\Omega^{\,-\alpha}$  outside  a larger ball.
Then $c_\Omega\,d_\Omega^{\;-2}\,|\varphi_\alpha|^2$ is integrable if $\alpha>(d+\delta'-2)/2$.
So choosing $\alpha=(d+\delta'-2+\varepsilon)/2$ and proceeding as in the local approximation one deduces that  $a\leq a(\delta')$.
Here the property $\lim_{s\to\infty}c(s)\,s^{-\delta'}=1$ is crucial.
\hfill$\Box$

\medskip

Now one can apply Theorem~\ref{tloc1.2} to obtain the Rellich inequality. 
It suffices to compute $\nu$ and verify that $\nu<a_1$.
There are two distinct cases.

\begin{obs}\label{oloc5.2} Assume $\delta+\delta'\leq4$.
If $d+2\,(\delta\wedge\delta')-4>0$ then the Rellich inequality
\begin{equation}
(H\varphi, H\varphi)\geq a_2\, (c_\Omega \,d_\Omega^{\;-2}\varphi,c_\Omega \,d_\Omega^{\;-2}\varphi)
\label{eloc5.2}
\end{equation}
is valid for all $\varphi\in D(H)$ with $a_2=d^{\,2}(d+2\,(\delta\wedge\delta')-4)^2/16$.
The value of $a_2$ is optimal.
\end{obs}
\proof\
First consider the case $\delta,\delta'\in[0,2]$.
Then
\[
\nu=\sup\{ |1-t/2|^2: \delta\wedge\delta'\leq t\leq \delta\vee\delta'\}=(1-(\delta\wedge\delta')/2)^2
\]
and $\nu<a_1$ if and only if $1-(\delta\wedge\delta')/2<(d+(\delta\wedge\delta')-2)/2$ or, equivalently, 
$d+2(\delta\wedge\delta')-4>0$.
But this implies the condition  $d+(\delta\wedge\delta')-2>0$ necessary for the  Hardy inequality (\ref{eloc5.0}). 
Therefore  one deduces from Theorem~\ref{tloc1.2} that  the Rellich inequality (\ref{eloc1.7}) is valid  with  constant $a_2=(a_1-\nu)^2$ 
which is easily calculated to be the value given in
the observation.

Secondly, assume $\delta\in[0,2]$ but $\delta'\geq2$.
Then $0\leq\delta'/2-1\leq 1-\delta/2$ since $\delta+\delta'\leq4$.
Therefore 
\[
\nu=\sup\{ |1-t/2|^2: \delta\leq t\leq \delta'\}=|1-\delta/2|^2=|1-(\delta\wedge\delta')/2|^2
\;.
\]
Hence the Rellich inequality (\ref{eloc5.2}) is again valid with the same value of $a_2$.

Thirdly, if $\delta'\in[0,2]$ but $\delta\geq2$ one reaches the same conclusion by interchanging
$\delta$ and $\delta'$ in the last argument.

Therefore the observation is established for all $\delta,\delta'\geq0$ with $\delta+\delta'\leq4$.

The optimality of $a_2$ follows by a  reasoning similar to the Hardy case. One again needs  separate
arguments at the origin and at infinity.
\hfill$\Box$

\begin{obs}\label{oloc5.3} Assume $\delta+\delta'\geq4$.
If $d-|\delta-\delta'|>0$ then the Rellich inequality $(\ref{eloc5.2})$
is valid  with $a_2=(d-|\delta-\delta'|)^2(d+\delta+\delta'-4)^2/16$.
\end{obs}
\proof\
First  assume $\delta,\delta'\geq2$.
Then $\nu=(1-(\delta\vee\delta')/2)^2$.
Therefore $a_1> \nu$ if and only if $d>(\delta\vee\delta')-(\delta\wedge\delta')=|\delta-\delta'|$.
This condition also  implies $d+(\delta\wedge\delta')\geq (\delta\vee\delta')\geq2$. 
Therefore the Hardy inequality (\ref{eloc5.0}) is valid.
Moreover, $a_1-\nu=(d+\delta+\delta'-4)(d-|\delta-\delta'|)/4$
and   one deduces from Theorem~\ref{eloc1.2} that the  Rellich inequality (\ref{eloc5.2}) is valid  with
$a_2=(a_1-\nu)^2$ whenever $d>|\delta-\delta'|$.

Secondly, assume $\delta\leq2$ and $\delta'\geq2$.
Since $\delta+\delta'\geq4$ it follows that $1-\delta/2\leq \delta'/2-1$.
Therefore $\nu=(1-\delta'/2)^2=(1-(\delta\vee\delta')/2)^2$.
Then the observation follows again.

The final case  $\delta\geq2$ and $\delta'\leq2$ now
follows from the second case by interchanging $\delta$ and $\delta'$.
\hfill$\Box$
\end{exam}

There is one question left over in this discussion of Example~\ref{exloc4.1},
the optimality of   the value of $a_2$ in Observation~\ref{oloc5.3}.
It does follow from the argument in the proof of Observation~\ref{oloc5.2}
that the value in the case $\delta+\delta'\geq 4$ is less than or equal
to the value in the case $\delta+\delta'\geq 4$.
Moreover the two values are equal if and only if  $\delta=\delta'$ or $\delta+\delta'=4$.
This follows by noting that 
\[
d(d+2\,(\delta\wedge\delta')-4)-(d-|\delta-\delta'|)(d+\delta+\delta'-4)
=|\delta-\delta'|\,(\delta+\delta'-4)
\;.
\]
Therefore the optimal value in the case $\delta+\delta'\geq4$ is generally strictly
smaller than the value $d^{\,2}(d+\delta\wedge\delta'-4)^2/16$.

\medskip

As a second  illustration of the foregoing techniques we consider a general class of operators of Grushin type
which are  related to the classic situation described in Example~\ref{exloc4.1}.
These operators  differ somewhat from the standard Grushin operators.
Many of their properties, e.g.\ Gausian kernel bounds, Poincar\'e inequalities, etc.,  were previously established in \cite{RSi2} \cite{RSi2a}  \cite{RSi6a}.
Although the Grushin  operators are not directly covered by Theorem~\ref{tloc1.2}  similar conclusions can be drawn by a slight modification
of the proof of Theorem~\ref{tloc1.1}.

\begin{exam}\label{exloc4.3}
Let $\Omega=(\Ri^{d_1}\backslash\{0\})\times\Ri^{d_2}$ and set $x=(x_1,x_2)$ with $x_1\in\Ri^{d_1}$ and $x_2\in \Ri^{d_2}$.
Then $\partial\Omega=\{x\!=\!(0, x_2):x_2\in \Ri^{d_2}\}$ and $d_\Omega(x)=|x_1|$.
Next define the Dirichlet  form $h$ on $L_2(\Omega)=L_2(\Ri^{d_1}\backslash\{0\})\otimes L_2(\Ri^{d_2})$  as the closure of the form
\begin{equation}
\varphi\in C_c^\infty(\Omega)\mapsto h(\varphi)=(\nabla_{\!x_1}\varphi,c_\Omega\,\nabla_{\!x_1}\varphi) +(\nabla_{\!x_2}\varphi, b\,\nabla_{\!x_2}\varphi)
\label{eloc5.50}
\end{equation}
where $c_\Omega=c\circ d_\Omega$ with $c$  again the function defined in Theorem~\ref{tloc1.2}
and $b$  the operator of multiplication by a positive bounded function of the first variable $x_1$.
Thus the coefficient $c_\Omega$ of the first form on the right, $h^{(1)}$,   and the coefficient $b$ of the second form, $h^{(2)}$,
are both  independent of $x_2$.

The forms $h^{(1)}$ and $h^{(2)}$ are both closable on $L_2(\Omega)$ and their closures are Dirichlet forms.
The submarkovian operator $H_1$ associated with $h^{(1)}$  is the tensor product of an operator $\widetilde H_1$ which acts on the first component  $L_2(\Ri^{d_1}\backslash\{0\})$ of the tensor product space and the identity operator 
$\one_2$ on the second component $L_2(\Ri^{d_2})$.
But $\widetilde H_1$  can be identified as the operator analyzed in Example~\ref{exloc4.1}.
Therefore it satisfies the Hardy inequality (\ref{eloc5.0}) on $L_2(\Ri^{d_1}\backslash\{0\})$, with $d$ replaced by $d_1$.
Now since $h\geq h^{(1)}$ the form $h$ satisfies the corresponding Hardy inequality on $L_2(\Omega)$. 
Explicitly one has the following.

\begin{obs}\label{oloc5.4}
If $d_1+(\delta\wedge\delta')-2>0$ then
\begin{equation}
h(\varphi)\geq a_1 \,(c_\Omega^{\,1/2}d_\Omega^{\;-1}\varphi, c_\Omega^{\,1/2}d_\Omega^{\;-1}\varphi)
\label{eloc5.51}
\end{equation}
 for all $\varphi\in D(h)$ with $a_1=(d_1+(\delta\wedge\delta')-2)^2/4$.
The value of $a_1$ is optimal.
\end{obs}
\proof\
It only remains to prove that the constant $a_1$ is optimal.
But if $\tilde a$ is the optimal constant then clearly $\tilde a\geq a_1$ and it suffices to prove the converse bound.

The optimal value $\tilde a$ is given by
\begin{eqnarray*}
\tilde a&=&\inf \{h(\varphi)/\|c_\Omega^{\,1/2}d_\Omega^{\;-1}\varphi\|_2^2:\varphi\in D(h)\}\\[5pt]
&\leq &\inf \{h(\psi\,\chi)/(\|c_\Omega^{\,1/2}d_\Omega^{\;-1}\psi\|_2^2\,\|\chi\|_2^2):\psi\in D(h^{(1)}), \chi\in C_c^\infty(\Ri^{d_2})\}
\end{eqnarray*}
where we have slightly abused notation by not distinguishing between the $L_2$-norms on the two components in the tensor product space.
It follows, however, from the product structure  that 
\[
h(\psi\,\chi)=h^{(1)}(\psi)\,\|\chi\|_2^2+(b\,\psi,\psi)\,\|\nabla_{\!x_2}\chi\|_2^2
\;.
\]
Next replace $\chi$ by $\chi_\lambda: \chi_\lambda(x_2)=\lambda^{d_2/2}\chi(\lambda x_2)$.
Since  $\|\chi_\lambda\|_2=\|\chi\|_2$ and $\|\nabla_{\!x_2}\chi_\lambda\|_2=\lambda\,\|\nabla_{\!x_2}\chi\|_2$ it follows 
immediately that 
\[
\lim_{\lambda\to0}h(\psi\,\chi_\lambda)/(\|c_\Omega^{\,1/2}d_\Omega^{\;-1}\psi\|_2^2\,\|\chi_\lambda\|_2^2)
=h^{(1)}(\psi)/\|c_\Omega^{\,1/2}d_\Omega^{\;-1}\psi\|_2^2=a_1
\]
where the last identification follows from Example~\ref{exloc4.1}.
Therefore $\tilde a=a_1$.
\hfill$\Box$.
\medskip

Next we argue that  the submarkovian operator $H$ corresponding to the Grushin form $h$ satisfies 
a Rellich inequality.  
Theorem~\ref{tloc1.2} is not directly applicable as the Grushin form has the second component $h^{(2)}$.
But in fact the Rellich inequality is independent of $h^{2)}$. 
This is somewhat surprising but can be understood by revisiting the proof of Theorem~\ref{tloc1.1}.
First we state the result.
There are again  two distinct regimes.

\begin{obs}\label{oloc5.5} 
Let $H$ be the submarkovian operator on $L_2(\Omega)$ corresponding to the Grushin form $(\ref{eloc5.50})$.

 Assume $\delta+\delta'\leq4$.
If $d_1+2\,(\delta\wedge\delta')-4>0$ then the Rellich inequality
\begin{equation}
(H\varphi, H\varphi)\geq a_2\, (c_\Omega \,d_\Omega^{\;-2}\varphi,c_\Omega \,d_\Omega^{\;-2}\varphi)
\label{eloc5.222}
\end{equation}
is valid for all $\varphi\in D(H)$ with $a_2=d_1^{\,2}(d_1+2\,(\delta\wedge\delta')-4)^2/16$.
The value of $a_2$ is optimal.

Alternatively assume $\delta+\delta'\geq4$.
If $d-|\delta-\delta'|>0$ then the Rellich inequality $(\ref{eloc5.222})$
is valid  with $a_2=(d_1-|\delta-\delta'|)^2(d_1+\delta+\delta'-4)^2/16$.
\end{obs}
\proof\
The Rellich inequality follows from the Hardy inequality of Observation~\ref{oloc5.4}
by a modification of the proof of Theorem~\ref{tloc1.1} with $\ce=h=h^{(1)}+h^{(2)}$.
The idea is to show that the main estimates of the proof are all independent of $h^{(2)}$.

First $h$ satisfies the $\eta$-Hardy inequality with $\eta=c_\Omega^{\,1/2}d_\Omega^{\;-1}$ by Observation~\ref{oloc5.4}
and $\eta$ is a function of  $x_1$.
Secondly, the form $h$ has a {\it carr\'e du champ} $\Gamma$ given by 
\[
\Gamma(\varphi)=c_\Omega\,|\nabla_{\!x_1}\varphi|^2+b\,|\nabla_{\!x_2}\varphi|^2
\;.
\]
Then since $\eta$ is independent of $x_2$ one has $\Gamma(\eta)=c_\Omega\,|\nabla_{\!x_1}\eta|^2$.
Thirdly it follows from the proof of Theorem~\ref{tloc1.2} that there exists an approximate identity $\rho_\alpha$
satisfying Condition~\ref{tloc1.1-3} of Theorem~\ref{tloc1.1} on $L_2(\Ri^{d_1}\backslash\{0\})$.
Therefore we can construct the bounded approximants $\eta_{\alpha, \beta}$ as in the proof of the latter theorem
and these remain functions of $x_1$.

The key identity (\ref{eloc201}) in the proof of Theorem~\ref{tloc1.1} now takes the form
\[
(H\varphi, \eta_{\alpha,\beta}^{\;2}\,\varphi_\varepsilon)=h(\eta_{\alpha,\beta}\,\varphi_\varepsilon)-h_{\varphi_\varepsilon^2}(\eta_{\alpha,\beta})
+h(\varphi-\varphi_\varepsilon, \eta_{\alpha,\beta}^{\;2}\,\varphi_\varepsilon)
\]
with $\varphi_\varepsilon$ again the bounded approximate to $\varphi\in D(H)$.
But $h\geq h^{(1)}$ and 
\[
h_{\varphi_\varepsilon^2}(\eta_{\alpha,\beta})=(\varphi, \Gamma(\eta_{\alpha,\beta})\varphi)=h^{(1)}_{\varphi_\varepsilon^2}(\eta_{\alpha,\beta})
\]
because $ \Gamma(\eta_{\alpha,\beta})$ is independent of $x_2$.
Therefore one obtains the estimate
\begin{eqnarray*}
\|H\varphi\|_2\,\| \eta_{\alpha,\beta}^{\;2}\,\varphi_\varepsilon\|_2&\geq& h^{(1)}(\eta_{\alpha,\beta}\,\varphi_\varepsilon)-h^{(1)}_{\varphi_\varepsilon^2}(\eta_{\alpha,\beta})
+h(\varphi-\varphi_\varepsilon, \eta_{\alpha,\beta}^{\;2}\,\varphi_\varepsilon)\\[5pt]
&\geq&h^{(1)}(\eta_{\alpha,\beta}\,\varphi_\varepsilon)-h^{(1)}_{\varphi_\varepsilon^2}(\eta_{\alpha,\beta})
-h^{(1)}(\varphi-\varphi_\varepsilon)^{1/2}\, h^{(1)}(\eta_{\alpha,\beta}^{\;2}\,\varphi_\varepsilon)^{1/2}\\[5pt]
&&\hspace{2.6cm}{}-h^{(2)}(\varphi-\varphi_\varepsilon)^{1/2} \,h^{(2)}(\eta_{\alpha,\beta}^{\;2}\,\varphi_\varepsilon)^{1/2}
\end{eqnarray*}
and the only dependence on $h^{(2)}$ is in the last term which converges to zero as $\varepsilon\to0$.
This last point depends on the uniform bound  
$h^{(2)}(\eta_{\alpha,\beta}^{\;2}\,\varphi_\varepsilon)\leq \|\eta_\beta\|_\infty^4\,
 h^{(2)}(\varphi_\varepsilon)\leq \|\eta_\beta\|_\infty^4 \,h^{(2)}(\varphi)$.
 Therefore one can now repeat the proof of Theorem~\ref{tloc1.1} following the identity (\ref{eloc201}) to obtain a Rellich inequality 
 which is totally independent of $h^{(2)}$.
The Rellich inequality is determined by $h^{(1)}$.
 The end result is (\ref{eloc5.222}) with $a_2=(a_1-\nu)^2$ where $a_1$ is the Hardy constant given in Observation~\ref{oloc5.4}
 and $\nu$ is again given by $\sup\{|1-t/2|^2: \delta\wedge\delta'\leq t\leq \delta\vee\delta'\}$.
 The calculation is a repetition of  that given in Example~\ref{exloc4.1}.
 The only difference is that $h^{(1)}$ is now a form on the first component of the tensor product space $L_2(\Ri^{d_1}\backslash\{0\})\otimes L_2(\Ri^{d_2})$ but this makes no essential difference.
 
 The only remaining point to establish is the optimality of $a_2$ in the case that $\delta+\delta'\leq 4$.
 But this follows by a slight generalization of the argument used to prove the optimality of $a_1$ in Observation~\ref{oloc5.4}.
 Now, however, one uses the tensor product structure to note that 
 \[
\|H(\psi\,\chi)\|_2^2=\|H_1\psi\|_2^2\,\|\chi\|_2^2+2\,(b\,\psi,H_1\psi)\,\|\nabla_{\!x_2}\chi\|_2^2+\|b\,\psi\|_2^2\,\|\Delta_{x_2}\chi\|_2^2
\]
with $\psi\in L_2(\Ri^{d_1}\backslash\{0\})$ and $\chi\in L_2(\Ri^{d_2})$.
 \hfill$\Box$
 \end{exam}

\section*{Acknowledgements}
The author is indebted to Juha Lehrb\"ack for a continuing, informative,  correspondence on the Hardy and Rellich inequalities
together with several specific  suggestions about the current material.
His advice has been very helpful.


\end{document}